\documentclass[a4paper,oneside,11pt]{article}

\newtheorem{dfn}{Definition}[section]
\newtheorem{tw}[dfn]{Theorem}
\newtheorem{prop}[dfn]{Proposition}
\newtheorem{rem}[dfn]{Remark}

\newtheorem{cor}[dfn]{Corollary}

\usepackage{anysize}
\usepackage{verbatim}
\usepackage{array}
\usepackage{enumerate}
\usepackage{amssymb} % symbole AMS
\usepackage{amsmath}
\usepackage{fancyhdr}
\usepackage{euscript}
\usepackage{latexsym}

\author{Micha\l \ Barski \\ \small  Faculty of Mathematics, Cardinal Stefan Wyszy\'nski University in Warsaw, Poland\\
\small Faculty of Mathematics and Computer Science, University of Leipzig, Germany\\ \small{\it m.baran@uksw.edu.pl} \bigskip \\
\\
Jerzy Zabczyk
\\ \small Institute of Mathematics, Polish Academy of
Sciences,
     Warsaw,  Poland\\ \small{\it zabczyk@impan.pl}}

\title{\bf Heath-Jarrow-Morton-Musiela  equation   with linear volatility\footnote{Supported by The Polish MNiSW grant
NN201419039.}}
\begin{document}
\baselineskip=1.2\baselineskip \maketitle
\date
\begin{abstract}
The paper is concerned with the problem of existence of solutions
for the Heath-Jarrow-Morton-Musiela equation with linear volatility.
Necessary conditions and sufficient conditions for the existence of
weak  and strong solutions are provided. The key role is played by  logarithmic growth conditions of
the L\'evy   exponent of the noise process  introduced in \cite{BarZab}.
\end{abstract}

%\tableofcontents

\section{Introduction}

Let $P(t,T)$ denote a price at time $t\geq0$ of a bond paying $1$
unit of money to its holder at time $T\geq t$. The prices
$P(\cdot,T)$ are processes defined on a fixed filtered probability
space $(\Omega, \mathcal{F}_{t, t\geq0},P)$. The forward rate $f$ is
a random field defined by the formula
\begin{gather*}
P(t,T)= e^{-\int_{t}^{T}f(t,u)du}, \qquad 0\leq t\leq T\leq T^\ast.
\end{gather*}
The prices of all bonds traded on the market are thus determined by
the forward rate $f(t,T), 0\leq t\leq T<+\infty$ and thus the
starting point in the bond market description is specifying the
dynamics of $f$. In this paper we consider the following stochastic
differentials
\begin{gather}\label{rownanie na f}
df(t,T)=\alpha(t,T)dt+\sigma(t,T)dL(t),\qquad 0\leq t\leq T,
\end{gather}
where $L$ is a L\'evy process. The equation above can be viewed as a
system of infinitely many equations parameterized by $0\leq
T<+\infty$. The discounted bond prices $\hat{P}(t,T)$ are defined by
\begin{gather*}
\hat{P}(t,T):=e^{-\int_{0}^{t}r(s)ds}\cdot P(t,T),\qquad 0\leq t\leq
T<+\infty,
\end{gather*}
where $r(t):=f(t,t), t\geq 0$ is the short rate. If we extend the
domain of $f$ by putting $f(t,T)=f(T,T)$ for $t\geq T$ we obtain the
formula
\begin{gather*}
\hat{P}(t,T)=e^{-\int_{0}^{T}f(t,u)du}, \qquad 0\leq t\leq T<
+\infty.
\end{gather*}
The market is supposed to be arbitrage free, i.e. we assume that the
processes $\hat{P}(\cdot,T)$ are local martingales. This implies
that the coefficients $\alpha,\sigma$ in \eqref{rownanie na f}
satisfy the Heath-Jarrow-Morton condition, i.e. for each $T\geq0$
\begin{gather}\label{warunek HJM}
\int_{t}^{T}\alpha(t,u)du=J\left(\int_{t}^{T}\sigma(t,u)du\right),
\end{gather}
for almost all $t\geq0$, see \cite{bjork}, \cite{Eberlein},
\cite{JakubowskiZabczyk}. The function $J$ above is the L\'evy
exponent of $L$ defined by
\begin{gather*}\label{J}
\mathbf{E}(e^{-zL(t)})=e^{tJ(z)},\qquad t\in[0,T^\ast], \ z\in
\mathbb{R},
\end{gather*}
where
\begin{gather}\label{Laplace transform}
J(z)=-az+\frac{1}{2}qz^2+\int_{\mathbb{R}}(e^{-zy}-1+zy\mathbf{1}_{(-1,
1)}(y)) \ \nu(dy),\qquad z\in\mathbb{R},
\end{gather}
with $ a\in \mathbb{R}$, $q\geq 0$ and the {\it L\'evy measure}
$\nu$ satisfying the integrability condition
\begin{equation}\label{war na miare Levyego}
\int_{\mathbb{R}}(y^2\wedge \ 1) \ \nu(dy)<\infty.
\end{equation}
Moreover, $J(z)$ is a finite number if and only if  $\int_{|y| \geq
1}(e^{-zy}) \ \nu(dy)<\infty$. As $J$ is differentiable,
\eqref{warunek HJM} can be written as
\begin{gather*}
\alpha(t,T)=J^\prime\left(\int_{t}^{T}\sigma(t,u)du\right)\sigma(t,T),
\quad 0\leq t\leq T<+\infty,
\end{gather*}
which means that the drift is fully determined by the volatility
process. As a consequence  \eqref{rownanie na f} reads as
\begin{gather}\label{rownanie na f przy HJM}
f(t,T)=f(0,T)+\int_{0}^{t}J^{\prime}\left(\int_{s}^{T}\sigma(s,u)du\right)\sigma(s,T)ds+\int_{0}^{t}\sigma(s,T)dL(s),
\quad 0\leq t\leq T<+\infty.
\end{gather}
The arguments of $f$ are the running time $t$ and the maturity date
$T$. Alternative description of the forward rate is provided by the
Musiela parametrization
\begin{gather*}
r(t,x):=f(t,t+x), \quad t\geq0, x\geq0,
\end{gather*}
involving $t$ as above and the time to maturity $x:=T-t$. The
initial curve will be denoted by  $r_0(x):=r(0,x), \ x\geq0$. It is
often more convenient to work with $r$ instead of $f$ because the
family of functions $\{r(t,\cdot)\}_t$ has a common domain
independent of $t$. Let $\{S_t, t\geq 0\}$ be the semigroup of
shifts, i.e. for any function $h$
\begin{gather*}
 S_t(h)(x):=h(t+x),\,\,\,t\geq 0,\,\,x\geq 0 .
\end{gather*}
To simplify the notation we set
$\tilde{\sigma}(t,x):=\sigma(t,t+x)$. Then, in virtue of
\eqref{rownanie na f przy HJM} we have {\small
\begin{align}\label{semigroup equation-wyprowadzenie}\nonumber
&r(t,x)=f(t,t+x)=f(0,t+x)+\int_{0}^{t}J^{\prime}\left(\int_{s}^{t+x}\sigma(s,u)du\right)\sigma(s,t+x)ds+\int_{0}^{t}\sigma(s,t+x)dL(s)\\[2ex]\nonumber
&=r(0,t+x)+\int_{0}^{t}J^{\prime}\left(\int_{0}^{x+t-s}\sigma(s,s+v)dv\right)\sigma(s,s+t-s+x)ds+\int_{0}^{t}\sigma(s,s+t-s+x)dL(s)\\[2ex]\nonumber
&=r_0(t+x)+\int_{0}^{t}J^{\prime}\left(\int_{0}^{x+t-s}\tilde{\sigma}(s,v)dv\right)\tilde{\sigma}(s,t-s+x)ds+\int_{0}^{t}\tilde{\sigma}(s,t-s+x)dL(s)\\[2ex]
&=S_t\Big(r_0(x)\Big)+\int_{0}^{t}S_{t-s}\left(J^{\prime}\Big(\int_{0}^{x}\tilde{\sigma}(s,v)dv\Big)\tilde{\sigma}(s,x)\right)ds+\int_{0}^{t}S_{t-s}\Big(\tilde{\sigma}(s,x)\Big)dL(s).
\end{align}}
If we assume that $\tilde{\sigma}$ is some function of $r$, i.e.
$\tilde{\sigma}(t,x)=(\tilde{\sigma} \circ r)(t,x)$ then $r$
satisfies
\begin{align}\label{semigroup equation}\nonumber
r(t,x)=S_t\Big(r_0(x)\Big)&+\int_{0}^{t}S_{t-s}\left(J^{\prime}\Big(\int_{0}^{x}(\tilde{\sigma}\circ r)(s,v)dv\Big)(\tilde{\sigma}\circ r)(s,x)\right)ds\\[2ex]
&+\int_{0}^{t}S_{t-s}\Big((\tilde{\sigma}\circ r)(s,x)\Big)dL(s),
\end{align}
and thus it is a weak solution of the semilinear, stochastic  equation
\begin{gather}\label{strong equation}
dr(t,x)=\left(A\Big(r(t,x)\Big)+J^{\prime}\Big(\int_{0}^{x}(\tilde{\sigma}\circ
r)(t,v)dv\Big)(\tilde{\sigma}\circ
r)(t,x)\right)dt+(\tilde{\sigma}\circ r)(t,x)dL(t),
\end{gather}
where $A$ stands for the generator of the semigroup $S_t$, i.e.
$A\big(h(x)\big)=\frac{d}{dx}h(x)$.

The problem of existence of solutions to \eqref{semigroup equation}
and \eqref{strong equation} is  under active  investigation in recent
time, see \cite{FilTap}, \cite{FilipovicTappeTeichmann},
\cite{Marinelli}, \cite{PeszatZabczyk}, where special assumptions
are imposed on $\tilde{\sigma}$ to obtain existence results for (\ref{strong equation}) or equivalently to  (\ref{semigroup equation}). In this paper we study the case of linear
volatility, i.e. it is assumed that
\begin{gather}\label{linear volatility}
\tilde{\sigma}(t,x)=\tilde{\lambda}(t,x)r(t-,x), \quad t\geq0, x\geq
0,
\end{gather}
where $\tilde{\lambda}(\cdot,\cdot)$ is a deterministic function
satisfying certain regularity conditions.

We show that  if
$J^{\prime}$ satisfies the logarithmic growth condition
\begin{gather}\label{logarytmiczny wzrost J^prim}
 \limsup_{z\rightarrow\infty} \ \left(\ln z-\bar{\lambda}T^\ast
 J^{\prime}\left(z\right)\right)=+\infty,\qquad 0<T^\ast<+\infty,
\end{gather}
then the   equation \eqref{semigroup equation} has solutions
in the weighted Hilbert spaces of square integrable functions or in the weighted Hilbert spaces of functions
with square integrable first derivative, see Theorem \ref{tw o
istnieniu}. It is also shown that if
\begin{gather}\label{drugi moment miary}
\int_{0}^{+\infty}y^2\nu(dy)<+\infty,
\end{gather}
then solutions are unique, see Theorem \ref{tw o jedynosci}. Moreover we prove that if
$\tilde{\sigma}(t,x)=\tilde{\lambda}(t)r(t-,x)$ and \eqref{drugi
moment miary} holds then the solution is strong in the space of
functions with square integrable first derivative, i.e. it solves
\eqref{strong equation}, see Theorem \ref{tw o silnym rozwiazaniu}.

On the other hand we  show that if $J^{\prime}$ grows  faster than a third power
of the logarithm, i.e.
\begin{gather}\label{ogr na Jprim}
J^{\prime}(z)\geq a(\ln z)^3+b, \qquad \forall z>0.
\end{gather}
for some $a>0$, $b\in\mathbb{R}$, then there is no non-exploding
weak solution on any finite time interval, see Theorem \ref{tw
o eksplozjach}.

In Section \ref{Existence and characteristics of the noise} we
give explicit conditions which imply \eqref{logarytmiczny wzrost J^prim}
or \eqref{ogr na Jprim}. They are formulated  in terms of the
parameters  of the noise process and  provide a precise description of the class of
L\'evy processes appropriate for  linear models.

The results are obtained via the random field approach. This enabled
us to relax some    assumptions required by the direct   SPDE approach.
 Let us also stress the   the
logarithmic growth condition \eqref{logarytmiczny wzrost J^prim}
admits in the equation ( \ref{strong equation}) coefficients which do not satisfy Lipschitz nor have
linear growth, so our results cover  non-standard
equations.

The logarithmic growth conditions \eqref{logarytmiczny wzrost
J^prim}, \eqref{ogr na Jprim} were  introduced in \cite{BarZab} and
examined in the space of bounded random fields on a finite domain.
In this paper we admit infinite domain, i.e. the solution is a
random field with unbounded second parameter and belongs to some
Hilbert space.

\section{Problem formulation}
Let us start with the description of the function
$\tilde{\lambda}:\mathbb{R}_{+}\times\mathbb{R}_{+}\longrightarrow
\mathbb{R}$ appearing in \eqref{linear volatility}. Denote by
$T^\ast$, where $0<T^\ast<+\infty$, a time horizon of the model,
i.e. $t\in[0,T^\ast]$. We assume that $\tilde{\lambda}(\cdot,x)$ is
continuous for each $x\geq0$ and that $\underline{\lambda}>0$,
$\bar{\lambda}<+\infty$, where
\begin{gather}\label{ograniczenia lambdy}
\underline{\lambda}=\underline{\lambda}(T^\ast):=\inf_{0\leq t\leq
T^\ast, x\geq0}\tilde{\lambda}(t,x), \qquad
\bar{\lambda}=\bar{\lambda}(T^\ast):=\sup_{0\leq t\leq T^\ast,
x\geq0}\tilde{\lambda}(t,x).
\end{gather}

The key step in our approach is to reduce the semigroup formulation
\eqref{semigroup equation} to the more  tractable operator form.
\begin{prop}\label{prop o postaci operatorowej dla r}
The field $r$ is a solution to \eqref{semigroup equation} if and
only if it satisfies
\begin{gather}\label{rownanie operatorwe na r}
r(t,x)=\tilde{a}(t,x)e^{\int_{0}^{t}J^{\prime}(\int_{0}^{t-s+x}\tilde{\lambda}(s,v)r(s,v)dv)\tilde{\lambda}(s,t-s+x)ds},
\quad t\geq0, x\geq0,
\end{gather}
where
\begin{align}\label{wzor na tilde a}\nonumber
 \tilde{a}(t,x):=&r_0(t+x)e^{\int_{0}^{t}\tilde{\lambda}(s,t-s+x)dL(s)-\frac{q^2}{2}\int_{0}^{t}\tilde{\lambda}^2(s,t-s+x)ds}\\[2ex]
&\cdot\prod_{0\leq s\leq t}\left(1+\tilde{\lambda}(s,t-s+x)\triangle
L(s)\right)e^{-\tilde{\lambda}(s,t-s+x)\triangle L(s)}.
\end{align}
\end{prop}
{\bf Proof:} The proof is based on the relation between $r$ and $f$.
It follows from \eqref{semigroup equation-wyprowadzenie},
\eqref{linear volatility} and \eqref{rownanie na f przy HJM} that
$r(t,x)=f(t,t+x), t\geq0, x\geq0$ is a solution of \eqref{semigroup
equation} if and only if $f(t,T), t\geq0, T\geq t$ satisfies {\small
\begin{gather}\label{rownanie na
f z linear volatility}
 f(t,T)=f(0,T)+\int_{0}^{t}J^{\prime}\left(\int_{s}^{T}\lambda(s,u)f(s-,u)du\right)\lambda(s,T)f(s-,T)ds+\int_{0}^{t}\lambda(s,T)
f(s-,T)dL(s),
\end{gather}}
with $\lambda(t,T):=\tilde{\lambda}(t,T-t)$. Now we see that $f$ is
a stochastic exponential and mimicking the proof of Proposition 2.1
in \cite{BarZab} we can show that $f$ solves \eqref{rownanie na f z
linear volatility} if and only if it satisfies
\begin{gather}\label{rownanie operatorowe na f}
 f(t,T) =
a(t,T)e^{\int_{0}^{t}J^\prime\left(\int_{s}^{T}\lambda(s,u)f(s,u)du\right)\lambda(s,T)ds},\quad
t\geq0,T\geq t,
\end{gather}
where
\begin{gather*}
a(t,T):=f_{0}(T)
 e^{\int_{0}^{t} \lambda (s, T) dL(s)-\frac{q^2}{2}\int_{0}^{t}\lambda^2(s,T)ds}\cdot\prod_{0\leq s\leq t}\Big(1+\lambda(s,T)\triangle L(s)\Big)e^{-\lambda(s,T)\triangle L(s)}.
\end{gather*}
If we put $T=t+x$ in \eqref{rownanie operatorowe na f} and check
that $\tilde{a}(t,x)=a(t,t+x)$ then we see that $f(t,T), t\geq 0,
T\geq t$ satisfies \eqref{rownanie operatorowe na f} if and only if
$r(t,x)=f(t,t+x), t\geq 0, x\geq0$ satisfies \eqref{rownanie
operatorwe na r}. \hfill$\square$

\subsection{Assumptions}

As forward rates are nonnegative, it is justified, in virtue
of \eqref{rownanie operatorwe na r}, \eqref{wzor na tilde a} and the
inequality $\tilde{\lambda}(t,x)<\bar{\lambda}$, to impose the
following standing assumptions:
\begin{enumerate}[]
\item (A1) \quad the initial curve $r_{0}$ is positive,
\item (A2) \quad the support of the L\'evy measure is contained in the interval
    $(-1/\bar{\lambda},+\infty)$.
\end{enumerate}
Moreover we will assume that
\begin{enumerate}[]
\item (A3) \quad
The random field $\Big\{\mid\int_{0}^{t}\tilde{\lambda}(s,t-s+x)dL(s)\mid:
(t,x)\Big\}$ is bounded on  $[0,T^\ast]\times[0,+\infty)$.
\end{enumerate}
Assumption $(A3)$ is satisfied for example for $\tilde{\lambda}$ of
the form
\begin{gather*}
\tilde{\lambda}(t,x)=\sum_{n=1}^{N}a_n(t)b_n(s+x),
\end{gather*}
where $\{a_n(\cdot)\}$ are continuous and $\{b_n(\cdot)\}$ are
bounded on $[0,+\infty)$.

In view of \eqref{rownanie operatorwe na r} it is clear that for our study of the equation (\ref{strong equation}) only the behavior of $J^{\prime}$ in
the interval $[0,+\infty)$ will be of interest. It is convenient (see, in particular Section \ref{Existence and characteristics of the
noise}) to decompose the function $J$ given by \eqref{Laplace
transform} to the from
\begin{gather}\label{J decomposition}
J(z)=-az+\frac{1}{2}qz^2+J_1(z)+J_2(z)+J_3(z),
\end{gather}
where
\begin{gather*}
J_1(z):=\int_{-1/{\bar {\lambda}}}^{0}(e^{-zy}-1+zy) \ \nu(dy),\quad
J_2(z):=\int_{0}^{1}(e^{-zy}-1+zy) \
\nu(dy)\\[2ex]
J_3(z):=\int_{1}^{\infty}(e^{-zy}-1) \ \nu(dy).
\end{gather*}
The functions $J_1,J_2,J_3$ are smooth on the interval
$(0,+\infty)$, see Lemma 8.1 and 8.2 in \cite{Rusinek}. Since in (\ref{strong equation}) we
need the first derivative of $J$ on the interval $[0,+\infty)$,
we additionally assume that $J^{\prime}(0)$ exists and is finite.
But, in view of the formula,
\begin{gather*}
J^{\prime}(0)=-a+J^{\prime}_3(0)=-a-\int_{1}^{+\infty}y\nu(dy),
\end{gather*}
this is equivalent to the assumption
$\int_{1}^{+\infty}y\nu(dy)<+\infty$. This, together with (A2)
gives the following standing assumption
\begin{gather*}
\text{(A4)}\qquad\qquad \int_{(-1/{\bar {\lambda}}, 1)}y^2 \nu(dy) +
\int_{1}^{\infty} y \nu (dy) <\infty.
\end{gather*}
It turns out, see Lemma 8.1 and 8.2 in \cite{Rusinek}, that under
$(A4)$ the functions  $J_1$, $J_2$, $J_3$ are smooth on the interval
$[0,+\infty)$ and
\begin{gather}\nonumber
J_1^{\prime}(z)=\int_{-1/{\bar {\lambda}}}^{0}y(1-e^{-zy})  \
\nu(dy),
\quad J_2^{\prime}(z)=\int_{0}^{1}y(1-e^{-zy})  \ \nu(dy)\\[2ex]\label{pierwsze pochodne}
J_3^{\prime}(z)=-\int_{1}^{\infty}ye^{-zy} \ \nu(dy),
\end{gather}
\begin{gather}\label{drugie pochodne}
J_1^{\prime\prime}(z)=\int_{-1/{\bar {\lambda}}}^{0}y^2e^{-zy} \
\nu(dy), \quad J_2^{\prime\prime}(z)=\int_{0}^{1}y^2e^{-zy} \
\nu(dy), \quad J_3^{\prime\prime}(z)=\int_{1}^{\infty}y^2e^{-zy} \
\nu(dy).
\end{gather}
It follows from \eqref{pierwsze pochodne} and \eqref{drugie
pochodne} that $J^{\prime}_1,J^{\prime}_2,J^{\prime}_3$ and thus
$J^{\prime}$ are increasing functions on $[0,+\infty)$.

\subsection{State spaces}
The forward rate $r$ is supposed to take values in the Hilbert
spaces defined below
\begin{align*}
L^{2,\gamma}_{+}&:=\{h:\mathbb{R}_{+}\longrightarrow\mathbb{R}_{+}:
\|h\|^{2}_{L^{2,\gamma}_{+}}:=\int_{0}^{+\infty}\mid
h(x)\mid^2e^{\gamma x}dx<+\infty\},\\[2ex]
H^{1,\gamma}_{+}&:=\{h:\mathbb{R}_{+}\longrightarrow\mathbb{R}_{+}:
\|h\|^{2}_{H^{1,\gamma}_{+}}:=\int_{0}^{+\infty}\left(\mid
h(x)\mid^2+\mid h^{\prime}(x)\mid^2\right)e^{\gamma x}dx<+\infty\}
\end{align*}
with $\gamma>0$. Thus we study the problem of existence of solution
to \eqref{rownanie operatorwe na r} such that
\begin{gather*}
r(t, \cdot)\in L^{2,\gamma}_{+} \quad \text{or} \quad r(t, \cdot)\in
H^{1,\gamma}_{+} \quad \text{for each} \ t\geq0.
\end{gather*}
The solution is called non-exploding in $L^{2,\gamma}_{+}$, resp.
$H^{1,\gamma}_{+}$ on the interval $[0,T^\ast]$ if
\begin{gather*}
\sup_{t\in[0,T^\ast]} \|r(t, \cdot)\|_{L^{2,\gamma}_{+}}<+\infty
\qquad \text{resp.} \qquad \sup_{t\in[0,T^\ast]}
\|r(t,\cdot)\|_{H^{1,\gamma}_{+}}<+\infty
\end{gather*}
with probability one. Recall that $T^\ast$ above stands for a time horizon of the model. By
$\mathbb{L}^{2,\gamma}_{+}=\mathbb{L}^{2,\gamma}_{+}(T^\ast)$,
$\mathbb{H}^{1,\gamma}_{+}=\mathbb{H}^{1,\gamma}_{+}(T^\ast)$ we
denote the space of functions
$h:[0,T^\ast]\times[0,+\infty)\longrightarrow\mathbb{R}_{+}$ with
finite norms
\begin{gather*}
\|h\|^2_{\mathbb{L}^{2,\gamma}_{+}}:=\sup_{t\in[0,T^\ast]}\|h(t,\cdot)\|^{2}_{L^{2,\gamma}_{+}},\qquad
\|h\|^2_{\mathbb{H}^{+}_{\gamma}}:=\sup_{t\in[0,T^\ast]}\|h(t,\cdot)\|^{2}_{H^{1,\gamma}_{+}}.
\end{gather*}
Thus $r$ is non-exploding in $L^{2,\gamma}_{+}$, resp.
$H^{1,\gamma}_{+}$ on the interval $[0,T^\ast]$ if and only if
\begin{gather*}
\|r\|^2_{\mathbb{L}^{2,\gamma}_{+}}<+\infty, \qquad \text{resp.}
\qquad \|r\|^2_{\mathbb{H}^{1,\gamma}_{+}}<+\infty,
\end{gather*}
with probability one.

Let us notice that if $h\in L^{2,\gamma}_{+}$ then
\begin{align}\label{calkowalnosc rozwiazania}\nonumber
\int_{0}^{+\infty}h(x)dx&=\int_{0}^{+\infty}h(x)e^{\frac{\gamma}{2}x}\cdot
e^{-\frac{\gamma}{2}x}dx\leq \left(\int_{0}^{+\infty}\mid
h(x)\mid^{2}e^{\gamma
x}dx\right)^{\frac{1}{2}}\left(\int_{0}^{+\infty}e^{-\gamma
x}dx\right)^\frac{1}{2}\\[2ex]
&\leq \frac{1}{\sqrt{\gamma}} \ \|h\|_{L^{2,\gamma}_{+}} <+\infty.
\end{align}
Thus the condition imposed on the forward rate
\begin{gather*}
\int_{0}^{+\infty}\mid r(t,x)\mid^2e^{\gamma x}dx<+\infty,
\end{gather*}
implies the non-degeneracy of bonds' prices at time $t$, i.e.
\begin{gather*}
P(t,T)=e^{-\int_{0}^{T-t}r(t,v)dv}>\varepsilon, \qquad T\geq t,
\end{gather*}
where $\varepsilon=\varepsilon(\omega,t)>0$. Consequently, if
$r\in\mathbb{L}^{2,\gamma}_{+}$ then the family of prices $\{P(t,T);
\ t\in[0,T^\ast], T\geq t\}$ is separated from zero uniformly in
$t$. The requirement
\begin{gather*}
\int_{0}^{+\infty}\mid r^{\prime}(t,x)\mid^2e^{\gamma x}dx<+\infty
\end{gather*}
is justified by the observations that the forward rates are getting
flat for large maturities. It is clear that if $h\in
H^{1,\gamma}_{+}$ then it is bounded. Indeed, the following
estimation holds
\begin{align}\nonumber
h(x)&=h(0)+\int_{0}^{x}h^{\prime}(y)dy\leq
h(0)+\int_{0}^{+\infty}h^{\prime}(y)e^{\frac{\gamma}{2}y}\cdot
e^{-\frac{\gamma}{2}y}dy\\[2ex]\nonumber
&\leq h(0)+\left(\int_{0}^{+\infty}\mid
h^{\prime}(y)\mid^{2}e^{\gamma
y}dy\right)^{\frac{1}{2}}\left(\int_{0}^{+\infty}e^{-\gamma
y}dy\right)^\frac{1}{2}\\[2ex]\label{ograniczonosc rozwiazania}
&\leq h(0)+\frac{1}{\sqrt{\gamma}} \ \|h\|_{H^{1,\gamma}_{+}},
\qquad \forall x\geq0.
\end{align}

To conclude, under $(A1) - (A4)$ we are searching for solutions of
\eqref{rownanie operatorwe na r} in the class of random fields
satisfying
\begin{gather*}
r(\cdot,x) \ \text{is adapted and c\`adl\`ag on} \ [0,T^\ast] \
\text{for each} \ x\geq0,\\[2ex]
P(r\in \mathbb{L}^{2,\gamma}_{+})=1, \quad \text{resp.} \quad P(r\in
\mathbb{H}^{1,\gamma}_{+})=1.
\end{gather*}

\section{Existence and uniqueness results}
Let $\mathcal{K}$ denote the operator, acting on functions of two
variables, defined by
\begin{gather}
\mathcal{K}h(t,x)=\tilde{a}(t,x)e^{\int_{0}^{t}J^{\prime}\left(\int_{0}^{t-s+x}\tilde{\lambda}(s,v)h(s,v)dv\right)\tilde{\lambda}(s,t+x)ds},
\quad t\geq0 , x\geq0,
\end{gather}
Then the equation \eqref{rownanie operatorwe na r} can be written in
the form $r=\mathcal{K}r$. The problem of existence of solutions
will be examined via properties of the iterative sequence of random
fields
\begin{gather}\label{schemat iteracyjny}
h_0\equiv0, \qquad h_{n+1}:=\mathcal{K}h_n, \qquad n=1,2,... \ .
\end{gather}
Let us write
$\tilde{a}$ in the form $\tilde{a}(t,x)=r_0(t+x)\tilde{b}(t,x)$,
where
\begin{align}\label{wzor na tilde b}\nonumber
\tilde{b}(t,x):=&e^{\int_{0}^{t}\tilde{\lambda}(s,t-s+x)dL(s)-
\frac{q^2}{2}\int_{0}^{t}\tilde{\lambda}^2(s,t-s+x)ds}\\[2ex]
&\cdot\prod_{0\leq s\leq t}\left(1+\tilde{\lambda}(s,t-s+x)\triangle
L(s)\right)e^{-\tilde{\lambda}(s,t-s+x)\triangle L(s)}.
\end{align}
It can be shown in the similar way as in the Proposition 2.3 in
\cite{BarZab} that under $(A1), (A2), (A3) $ the field $\tilde{b}$ is bounded, i.e.
\begin{gather}\label{ograniczonosc tilde b}
\sup_{t\in[0,T^\ast], x\geq0} \tilde{b}(t,x)<\bar{b},
\end{gather}
where $\bar{b}=\bar{b}(\omega)>0$.
It can be shown by induction that if $r_0\in L^{2,\gamma}_{+}$ then
$h_n\in\mathbb{L}^{2,\gamma}_{+}$ for each $n$. Indeed, if
$h_{n}\in\mathbb{L}^{2,\gamma}_{+}$ then, in view of
\eqref{calkowalnosc rozwiazania} and \eqref{ograniczonosc tilde b},
we have
\begin{align*}
h_{n+1}(t,x)&\leq
r_0(t+x)\ \bar{b}\ e^{\bar{\lambda}\int_{0}^{t}\mid J^{\prime}(\int_{0}^{t-s+x}\tilde{\lambda}(s,v)h_n(s,v)dv)\mid ds}\\[2ex]
&\leq r_0(t+x) \ \bar{b} \ e^{\bar{\lambda}T^\ast\big|
J^{\prime}(\frac{\bar{\lambda}}{\sqrt{\gamma}}\parallel
h_n\parallel_{\mathbb{L}^{2,\gamma}_{+}})\big|},
\end{align*}
and thus $h_{n+1}\in\mathbb{L}^{2,\gamma}_{+}$. It follows from
$(A1)$, the assumption $\underline{\lambda}>0$ and the fact that
$J^{\prime}$ is increasing that the sequence $\{h_n\}$ is
monotonically increasing and thus there exists
$\bar{h}:[0,T^\ast]\times[0,+\infty)\longrightarrow \mathbb{R}_{+}$
such that
\begin{gather}\label{granica w schemacie iteracyjnym}
\lim_{n\rightarrow +\infty}h_n(t,x)=\bar{h}(t,x), \qquad 0\leq t\leq
T^\ast, x\geq0.
\end{gather}
Passing to the limit in \eqref{schemat iteracyjny}, by the monotone
convergence, we obtain
\begin{gather*}
\bar{h}(t,x)=\mathcal{K}h(t,x), \qquad 0\leq t\leq T^\ast, x\geq0.
\end{gather*}
It turns out that properties of the field $\bar{h}$ strictly depend
on the growth of the function $J^{\prime}$. In Section \ref{section
Existence} we show that if \eqref{logarytmiczny wzrost J^prim} holds
then $\bar{h}\in\mathbb{L}^{2,\gamma}_{+}$. Additional assumptions
guarantee that $\bar{h}\in\mathbb{H}^{1,\gamma}_{+}$ and that the
solution is unique. In Section \ref{section Explosions} it is shown
that if $J^{\prime}$ satisfies \eqref{ogr na Jprim} then $\bar{h}$
with positive probability is not in
$\mathbb{L}^{2,\gamma}_{+}=\mathbb{L}^{2,\gamma}_{+}(T^\ast)$ for
any $T^\ast$ and consequently that any random field $r$ satisfying
$\mathcal{K}r=r$ is not in $\mathbb{L}^{2,\gamma}_{+}$.

\subsection{Existence of weak solutions}\label{section Existence}

We start with an auxiliary result.
\begin{prop}\label{prop pomocniczy}
Assume that $J^{\prime}$ satisfies \eqref{logarytmiczny wzrost
J^prim}. If $r_0\in L^{2,\gamma}_{+}$ then there exists a positive
constant $c_1$ such that if
\begin{gather*}
\|h\|_{\mathbb{L}^{2,\gamma}_{+}}\leq c_1
\end{gather*}
then
\begin{gather*}
\|\mathcal{K}h\|_{\mathbb{L}^{2,\gamma}_{+}}\leq c_1.
\end{gather*}
\end{prop}
{\bf Proof}: $a)$ By \eqref{calkowalnosc rozwiazania} and
\eqref{ograniczonosc tilde b}, for any $t\in[0,T^\ast]$, we have
\begin{align*}
\|\mathcal{K}h(t,\cdot)\|_{L^{2,\gamma}_{+}}^2&=\int_{0}^{+\infty}|r_0(t+x)\tilde{b}(t,x)|^2e^{2\int_{0}^{t}J^{\prime}(\int_{0}^{t-s+x}\tilde{\lambda}(s,v)h(s,v)dv)\tilde{\lambda}(s,t-s+x)ds}
e^{\gamma x}dx\\[2ex]
&\leq \bar{b}^2\int_{0}^{+\infty}|r_0(t+x)|^2e^{2
J^{\prime}\left(\frac{\bar{\lambda}}{\sqrt{\gamma}}\cdot\|h\|_{\mathbb{L}^{2,\gamma}_{+}}\right)\int_{0}^{t}\tilde{\lambda}(s,t-s+x)ds}e^{\gamma
x}dx\\[2ex]
&\leq \bar{b}^2\cdot \|r_0\|^2_{L^{2,\gamma}_{+}}\cdot\sup_{s\in[0,t],
x\geq0}
e^{2J^{\prime}\left(\frac{\bar{\lambda}}{\sqrt{\gamma}}\cdot\|h\|_{\mathbb{L}^{2,\gamma}_{+}}\right)\int_{0}^{t}\tilde{\lambda}(s,t-s+x)ds}.
\end{align*}
This implies
\begin{gather*}
\|\mathcal{K}h\|_{\mathbb{L}^{2,\gamma}_{+}}\leq \bar{b} \cdot
\|r_0\|_{L^{2,\gamma}_{+}} \cdot \sup_{t\in[0,T^\ast], s\in[0,t],
x\geq0}e^{
J^{\prime}\left(\frac{\bar{\lambda}}{\sqrt{\gamma}}\cdot\|h\|_{\mathbb{L}^{2,\gamma}_{+}}\right)\int_{0}^{t}\tilde{\lambda}(s,t-s+x)ds},
\end{gather*}
and thus it is enough to find constant $c_1$ such that
\begin{gather}\label{warunek na c}
\ln\left(\bar{b}\cdot\|r_0\|_{L^{2,\gamma}_{+}}\right)+
\sup_{t\in[0,T^\ast], s\in[0,t],
x\geq0}J^{\prime}\left(\frac{\bar{\lambda}c_1}{\sqrt{\gamma}}\right)\int_{0}^{t}\tilde{\lambda}(s,t-s+x)ds\leq
\ln c_1.
\end{gather}
If $J^{\prime}(z)\leq 0$ for each $z\geq0$ then we put
$c_1=\bar{b}\cdot\|r_0\|_{L^{2,\gamma}_{+}}$. If $J^{\prime}$ takes
positive values then it is enough to find large $c_1$ such that
\begin{gather*}
\ln\left(\bar{b}\cdot\|r_0\|_{L^{2,\gamma}_{+}}\right)\leq \ln c_1 -
\bar{\lambda}T^\ast
J^{\prime}\left(\frac{\bar{\lambda}c_1}{\sqrt{\gamma}}\right).
\end{gather*}
Existence of such $c_1$ is a consequence of \eqref{logarytmiczny
wzrost J^prim}.\hfill$\square$\\

For the next result we will need to impose additional assumption on
the regularity of $\tilde{\lambda}$, i.e. that
$\tilde{\lambda}(t,\cdot)$ and $\tilde{b}(t,\cdot)$ are
differentiable and
\begin{gather}\label{ograniczenie pochodnej tilde lambda}
\sup_{t\in[0,T^\ast],x\geq0}\mid\tilde{\lambda}^{\prime}_x(t,x)\mid<+\infty,\\[2ex]\label{ograniczenie pochodnej tilde b}
\sup_{t\in[0,T^\ast],x\geq0}\mid\tilde{b}^{\prime}_x(t,x)\mid<+\infty.
\end{gather}

\begin{tw}\label{tw o istnieniu}
Assume that conditions (A.1) to (A.4) and  $\eqref{logarytmiczny wzrost J^prim}$ hold.
\begin{enumerate}[a)]
\item If $r_0\in L^{2,\gamma}_{+}$ then there exists a solution to
\eqref{rownanie operatorwe na r} taking values in the space
$L^{2,\gamma}_{+}$.
\item Assume that \eqref{ograniczenie pochodnej tilde lambda} and \eqref{ograniczenie pochodnej tilde b} are satisfied.
If $r_0\in H^{1,\gamma}_{+}$ and \eqref{drugi moment miary} holds then there exists a solution to \eqref{rownanie operatorwe na r} taking values
in the space $\mathbb{H}^{1,\gamma}_{+}$. 
\end{enumerate}
\end{tw}
{\bf Proof}: The limit $\bar{h}(\cdot,x)$ is adapted for each
$x\geq0$ as a pointwise limit.\\
\noindent $(a)$  Let $c_1$ be a
constant given by Proposition \ref{prop pomocniczy}. Then the
sequence $\{h_n\}$ is bounded in $\mathbb{L}^{2,\gamma}_{+}$ and
thus by the Fatou lemma we have
\begin{gather*}
\sup_{t\in[0,T^\ast]}\int_{0}^{+\infty}\mid \bar{h}(t,x)\mid^2
e^{\gamma x}dx\leq \sup_{t\in[0,T^\ast]}\underset{n\rightarrow
+\infty}{\lim \inf}\int_{0}^{+\infty}\mid h_n(t,x)\mid^2 e^{\gamma
x}dx\leq c_1^2,
\end{gather*}
and hence $\bar{h}\in\mathbb{L}^{2,\gamma}_{+}$.

$b)$ We will show that the solution $\bar{h}$ belongs to
$\mathbb{H}^{1,\gamma}_{+}$. Differentiating the equation $\bar
h=\mathcal{K}h$ gives
\begin{gather*}
\bar{h}^{\prime}(t,x)=r_0^{\prime}(t+x)\tilde{b}(t,x)F_1(t,x)+r_0(t+x)\tilde{b}^{\prime}_x(t,x)F_1(t,x)+r_0(t+x)\tilde{b}(t,x)F_1(t,x)F_2(t,x),
\end{gather*}
where
\begin{align*}
F_1(t,x)&:=e^{\int_{0}^{t}J^{\prime}(\int_{0}^{t-s+x}\tilde{\lambda}(s,v)\bar{h}(s,v)dv)\tilde{\lambda}(s,t-s+x)ds},\\[2ex]
F_2(t,x)&:=\int_{0}^{t}J^{''}\left(\int_{0}^{t-s+x}\tilde{\lambda}(s,v)\bar{h}(s,v)dv\right)\tilde{\lambda}^2(s,t-s+x)\bar{h}(s,t-s+x)ds\\[2ex]
&+\int_{0}^{t}J^{\prime}\left(\int_{0}^{t-s+x}\tilde{\lambda}(s,v)h(s,v)dv\right)\tilde{\lambda}^{\prime}_{x}(s,t-s+x)ds.
\end{align*}
We will show that
\begin{gather*}
\sup_{t\in[0,T^\ast], x\geq 0}F_1(t,x)<+\infty, \quad
\sup_{t\in[0,T^\ast], x\geq 0}F_2(t,x)<+\infty.
\end{gather*}
Then, in view of \eqref{ograniczonosc tilde b}, \eqref{ograniczenie
pochodnej tilde b}, the assertion follows from the assumption
$r_0\in H^{1,\gamma}_{+}$.

\noindent
From the fact $\bar{h}\in\mathbb{L}^{2,\gamma}_{+}$ it follows that
\begin{gather*}
\sup_{t\in[0,T^\ast], x\geq 0}F_1(t,x)\leq e^{\big|
J^{\prime}\left(\frac{\bar{\lambda}}{\sqrt{\gamma}}\parallel
\bar{h}\parallel_{\mathbb{L}^{2,\gamma}_{+}}\right)\big|\bar{\lambda}T^\ast}<+\infty.
\end{gather*}
It can be shown, see Theorem \ref{tw o warunkach koniecznych}, that
if $L$ admits negative jumps or contains Wiener part then \eqref{ogr
na Jprim} holds and consequently \eqref{logarytmiczny wzrost J^prim}
dos not hold. Thus \eqref{logarytmiczny wzrost J^prim} implies that
$J^{''}$ reduces to the form
$J^{''}(z)=\int_{0}^{+\infty}y^2e^{-zy}\nu(dy)$ and $0\leq
J^{''}(0)<+\infty$ due to the assumption \eqref{drugi moment miary}.
Since $J^{''}$ is decreasing, the following estimation holds
\begin{align*}
\sup_{t\in[0,T^\ast], x\geq 0}F_2(t,x)\leq&
J^{''}(0)T^\ast\bar{\lambda}^2\sup_{t\in[0,T^\ast], x\geq
0}\int_{0}^{t}\bar{h}(s,t-s+x)ds\\[2ex]
&+T^\ast\big|
J^{\prime}\left(\frac{\bar{\lambda}}{\sqrt{\gamma}}\parallel
\bar{h}\parallel_{\mathbb{L}^{2,\gamma}_{+}}\right)\big|\cdot\sup_{t\in[0,T^\ast],
x\geq0}\tilde{\lambda}^{\prime}_x(t,x).
\end{align*}
In view of \eqref{ograniczenie pochodnej tilde lambda} it is enough
to show that $\bar{h}$ is bounded on $\{(t,x), t\in[0,T^\ast],
x\geq0\}$. Using the fact that $\bar{h}=\mathcal{K}\bar{h}$ and
\eqref{ograniczonosc rozwiazania} we obtain
\begin{gather*}
\sup_{t\in[0,T^\ast], x\geq
0}\bar{h}(t,x)\leq\sup_{x\geq0}r_0(x)\cdot\sup_{t\in[0,T^\ast],
x\geq 0}\tilde{b}(t,x)\cdot
e^{\big|J^{\prime}\left(\frac{1}{\sqrt{\gamma}}\parallel\bar{h}\parallel_{L^{2,\gamma}_{+}}\right)\big|\bar{\lambda}T^\ast}<+\infty.
\end{gather*}
\hfill $\square$

\begin{rem}
Let $w:\mathbb{R}_{+}\longrightarrow\mathbb{R}_{+}$ be such that
$\int_{0}^{+\infty}\frac{1}{w(x)}dx<+\infty$. It can be shown with
similar proofs that the condition \eqref{logarytmiczny wzrost
J^prim} implies existence of non-exploding solution of
\eqref{rownanie operatorwe na r} taking values in the spaces
\begin{align*}
L^{2+}_w&:=\{h:\mathbb{R}_{+}\longrightarrow\mathbb{R}_{+}:
\int_{0}^{+\infty}\mid h(x)\mid^2 w(x)dx<+\infty\},
\end{align*}
and if $r_0$ is bounded and \eqref{drugi moment miary} holds, then
also in the space
\begin{align*}
H^{+}_w&:=\{h:\mathbb{R}_{+}\longrightarrow\mathbb{R}_{+}:
\int_{0}^{+\infty}\left(\mid h(x)\mid^2+\mid
h^{\prime}(x)\mid^2\right)w(x)dx<+\infty\}.
\end{align*}
\end{rem}

\subsection{Existence of strong solutions}

Under additional conditions we can establish existence of strong solutions.
\begin{tw}\label{tw o silnym rozwiazaniu}
Assume that
\begin{gather}\label{lambda niezalezna od x}
\tilde{\lambda}(t,x)=\tilde{\lambda}(t) \quad for \ x\geq0, t\geq 0,
\end{gather}
$r_0\in H^{1,\gamma}_{+}$ and \eqref{drugi moment miary} holds. Then
the non-exploding solution given by Theorem \ref{tw o istnieniu} (b)
is a strong solution of \eqref{strong equation}.
\end{tw}
{\bf Proof:} Taking into account \eqref{lambda niezalezna od x} and
differentiating \eqref{rownanie operatorwe na r} provides 
{\small
\begin{align}\label{pochodna r}\nonumber
\frac{\partial}{\partial x} r(t,x)&=
e^{\int_{0}^{t}\tilde{\lambda}(s)dL_s-\frac{q^2}{2}\int_{0}^{t}\tilde{\lambda}^2(s)ds}\prod(1+\tilde{\lambda}(s)\triangle
L_s)e^{-\tilde{\lambda}(s)\triangle
L_s}\cdot\\[2ex]\nonumber
&\cdot\bigg(r^{\prime}_0(t+x)e^{\int_{0}^{t}J^{\prime}(\int_{0}^{t-s+x}\tilde{\lambda}(s)r(s,v)dv)\tilde{\lambda}(s)ds}+r_0(t+x)e^{\int_{0}^{t}J^{\prime}(\int_{0}^{t-s+x}\tilde{\lambda}(s)r(s,v)dv)\tilde{\lambda}(s)ds}\cdot\\[2ex]\nonumber
&\phantom{aaaaaaaaaaaaaaaaaaaaaa}
\cdot\int_{0}^{t}J^{''}\Big(\int_{0}^{t-s+x}\tilde{\lambda}(s)r(s,v)dv\Big)\cdot
\tilde{\lambda}^2(s)r(s,t-s+x)ds\bigg)\\[2ex]\nonumber
&=r(t,x)\frac{r^{\prime}_0(t+x)}{r_0(t+x)}+r(t,x)\int_{0}^{t}J^{''}\Big(\int_{0}^{t-s+x}\tilde{\lambda}(s)r(s,v)dv\Big)\cdot\tilde{\lambda}^2(s)
r(s,t-s+x)ds\\[2ex]
&=r(t,x)\bigg[ \
\frac{r^{\prime}_0(t+x)}{r_0(t+x)}+\int_{0}^{t}J^{''}\Big(\int_{0}^{t-s+x}\tilde{\lambda}(s)r(s,v)dv\Big)\cdot
\tilde{\lambda}^2(s)r(s,t-s+x)ds\bigg].
\end{align}}
For $Z_1$, $Z_2$ defined by
\begin{align*}
Z_1(t)&:=e^{\int_{0}^{t}\tilde{\lambda}(s)dL_s-\frac{q^2}{2}\int_{0}^{t}\tilde{\lambda}^2(s)ds}\prod(1+\tilde{\lambda}(s)\triangle
L_s)e^{-\tilde{\lambda}(s)\triangle
L_s},\\[2ex]
Z_2(t,x)&:=r_0(t+x)e^{\int_{0}^{t}J^{\prime}(\int_{0}^{t-s+x}\tilde{\lambda}(s)r(s,v)dv)\tilde{\lambda}(s)ds},
\end{align*}
we have SDE's of the form {\small
\begin{align*}
dZ_1(t)&=Z_1(t-)\tilde{\lambda}(t)dL(t)\\[2ex]
dZ_2(t,x)&=\bigg\{r^{\prime}_0(t+x)e^{\int_{0}^{t}J^{\prime}(\int_{0}^{t-s+x}\tilde{\lambda}(s)r(s,v)dv)\tilde{\lambda}(s)ds}+r_0(t+x)e^{\int_{0}^{t}J^{\prime}(\int_{0}^{t-s+x}\tilde{\lambda}(s)r(s,v)dv)\tilde{\lambda}(s)ds}\cdot\\[2ex]
&\cdot\Big[J^{\prime}\Big(\int_{0}^{x}\tilde{\lambda}(t)r(t,v)dv\Big)\tilde{\lambda}(t)+\int_{0}^{t}J^{''}\Big(\int_{0}^{t-s+x}\tilde{\lambda}(s)r(s,v)dv\Big)\tilde{\lambda}^2(s)r(s,t-s+x)ds\Big]\bigg\}dt\\[2ex]
&=\bigg\{\frac{r^{\prime}_0(t+x)}{r_0(t+x)}Z_2(t,x)+Z_2(t,x)\Big[J^{\prime}\Big(\int_{0}^{x}\tilde{\lambda}(t)r(t,v)dv\Big)\tilde{\lambda}(t)+\\[2ex]
&\phantom{aaaaaaaaaaaaaaaaaaaaaaaa}+\int_{0}^{t}J^{''}\Big(\int_{0}^{t-s+x}\tilde{\lambda}(s)r(s,v)dv\Big)\tilde{\lambda}^2(s)r(s,t-s+x)ds\Big]\bigg\}dt\\[2ex]
&=\bigg\{Z_2(t,x)\bigg[\frac{r^{\prime}_0(t+x)}{r_0(t+x)}
+J^{\prime}\Big(\int_{0}^{x}\tilde{\lambda}(t)r(t,v)dv\Big)\tilde{\lambda}(t)+\\[2ex]
&\phantom{aaaaaaaaaaaaaaaaaaaaaaaa}+\int_{0}^{t}J^{''}\Big(\int_{0}^{t-s+x}\tilde{\lambda}(s)r(s,v)dv\Big)\tilde{\lambda}^2(s)r(s,t-s+x)ds\bigg]\bigg\}dt.
\end{align*}}
Using the formulas above, we obtain SDE for $r(t,x)$: {\small
\begin{align*}
dr(t,x)&=d\Big(Z_1(t)Z_2(t,x)\Big)=Z_1(t)dZ_2(t,x)+Z_2(t,x)dZ_1(t)\\[2ex]
&=Z_1(t)Z_2(t,x)\bigg[\frac{r^{\prime}_0(t+x)}{r_0(t+x)}
+J^{\prime}\Big(\int_{0}^{x}\tilde{\lambda}(t)r(t,v)dv\Big)\tilde{\lambda}(t)+\\[2ex]
&\phantom{aaaaaaaaaaaaaaaaaaaaaa}+\int_{0}^{t}J^{''}\Big(\int_{0}^{t-s+x}\tilde{\lambda}(s)r(s,v)dv\Big)\tilde{\lambda}^2(s)r(s,t-s+x)ds\bigg]dt\\[2ex]
&+Z_2(t,x)Z_1(t-)\tilde{\lambda}(t)dL(t)\\[2ex]
&=r(t,x)\bigg[\frac{r^{\prime}_0(t+x)}{r_0(t+x)}
+\int_{0}^{t}J^{''}\Big(\int_{0}^{t-s+x}\tilde{\lambda}(s)r(s,v)dv\Big)\tilde{\lambda}^2(s)r(s,t-s+x)ds\bigg]dt\\[2ex]
&+r(t,x)J^{\prime}\Big(\int_{0}^{x}\tilde{\lambda}(t)r(t,v)dv\Big)\tilde{\lambda}(t)dt+r(t-,x)\tilde{\lambda}(t)dL(t)\\[2ex]
&\overset{by \eqref{pochodna r}}{=}\frac{\partial}{\partial x}
r(t,x)dt+J^{\prime}\Big(\int_{0}^{x}\tilde{\lambda}(t)r(t-,v)dv\Big)\tilde{\lambda}(t)r(t-,x)dt+r(t-,x)\tilde{\lambda}(t)dL(t),
\end{align*}}
which is \eqref{strong equation}. \hfill $\square$

\subsection{Uniqueness}

In the next part of this section we investigate the problem of
uniqueness of solution.
\begin{prop}\label{prop Gronwall}
Let $d:\mathcal{P}\longrightarrow\mathbb{R}_{+}$ be a bounded
function satisfying
\begin{gather}\label{wzor w Gronwallu}
d(t,x)\leq K \int_{0}^{t}\int_{0}^{t-s+x}d(s,v)dvds,
\end{gather}
where $K>0$. Then $d(t,x)=0$ for all $(t,x)\in
[0,T^\ast]\times[0,+\infty)$.
\end{prop}
{\bf Proof:} Let $d$ be bounded by $M>0$ on
$[0,T^\ast]\times[0,+\infty)$. Let us define a new function
\begin{gather*}
\bar{d}(u,w):=d(u,w-u); \qquad u\in[0,T^\ast], w\geq u.
\end{gather*}
It is clear that $d\equiv0$ on $[0,T^\ast]\times[0,+\infty)$ if and
only if $\bar{d}\equiv 0$ on the set $\{(u,w): u\in[0,T^\ast], w\geq
u\}$. Let us notice that \eqref{wzor w Gronwallu} implies that
\begin{align*}
\bar{d}(u,w)&=d(u,w-u)\leq K\int_{0}^{u}\int_{0}^{w-s}d(s,y)dy
ds\\[2ex]
&= K \int_{0}^{u}\int_{s}^{w}d(s,z-s)dz ds=K
\int_{0}^{u}\int_{s}^{w}\bar{d}(s,z)dz ds.
\end{align*}
Using this inequality we will show by induction that
\begin{gather}\label{indukcja 2}
\bar{d}(u,w)\leq MK^n\frac{(uw)^n}{(n!)^2}, \qquad n=0,1,2,... .
\end{gather}
Then letting $n\rightarrow0$ we have $\bar{d}(t,x)=0$. The formula
\eqref{indukcja 2} is valid for $n=0$. Assume that it is true for
$n$ and show for $n+1$.
\begin{align*}
\bar{d}(u,w)&\leq K\int_{0}^{u}\int_{s}^{w}MK^{n}\frac{(sz)^n}{(n!)^2}dzds=MK^{n+1}\frac{1}{(n!)^2}\int_{0}^{u}s^n(\int_{s}^{w}z^n dz)ds\\[2ex]
&=
MK^{n+1}\frac{1}{(n!)^2}\int_{0}^{u}s^n\left(\frac{w^{n+1}-s^{n+1}}{n+1}\right)ds\leq
MK^{n+1}\frac{1}{(n!)^2}\int_{0}^{u}s^n\frac{w^{n+1}}{n+1}ds\\[2ex]
&=MK^{n+1}\frac{1}{(n!)^2}\frac{u^{n+1}}{(n+1)}\frac{w^{n+1}}{(n+1)}=MK^{n+1}\frac{(uw)^{n+1}}{((n+1)!)^2}.
\end{align*}
\hfill $\square$

\begin{tw}\label{tw o jedynosci}
Assume that $r_0^{\ast}:=\sup_{x\geq0}r_0(x)<+\infty$ and
\eqref{drugi moment miary} holds. If, on the interval $[0,T^\ast]$,
there exists  a non-exploding solution of the equation
\eqref{rownanie operatorwe na r}  taking values in
$L^{2,\gamma}_{+}$ then it is unique.
\end{tw}
{\bf Proof:} Assume that $r_1, r_2\in \mathbb{L}^{2+}_{+}$ are two
solutions of the equation \eqref{rownanie operatorwe na r} and
define
\begin{gather*}
d(t,x):=\mid r_{1}(t,x)-r_2(t,x)\mid, \qquad 0\leq t\leq T^\ast,
x\geq0.
\end{gather*}
Denote  $B:=\sup_{t\in[0,T^\ast],x\geq0}\tilde{b}(t,x)$. By
\eqref{rownanie operatorwe na r} and \eqref{calkowalnosc
rozwiazania}, for any $(t,x)\in[0,T^\ast]\times[0,+\infty)$, we have
\begin{align*}
d(t,x)&\leq
r_0(t+x)\tilde{b}(t,x)\left[e^{\int_{0}^{t}J^{\prime}(\int_{0}^{t-s+x}\tilde{\lambda}(s,v)r_1(s,v)dv)\tilde{\lambda}(s,t-s+x)ds}
+e^{\int_{0}^{t}J^{\prime}(\int_{0}^{t-s+x}\tilde{\lambda}(s,v)r_2(s,v)dv)\tilde{\lambda}(s,t-s+x)ds}\right]\\[2ex]
&\leq r_0^{\ast}\cdot B\cdot\left[e^{\bar{\lambda}T^\ast
\big|J^{\prime}(\frac{\bar{\lambda}}{\sqrt{\gamma}}\|r_1\|_{\mathbb{L}^{2,\gamma}_{+}})\big|}
+e^{\bar{\lambda}T^\ast
J^{\prime}\big|(\frac{\bar{\lambda}}{\sqrt{\gamma}}\|r_2\|_{\mathbb{L}^{2,\gamma}_{+}})\big|}
\right]<+\infty,
\end{align*}
and thus $d$ is bounded on $[0,T^\ast]\times[0,+\infty)$. In view of
the inequality $\mid e^{x}-e^{y}\mid\leq e^{x\vee y}\mid x-y\mid; \
x,y\geq0$ and the fact that $J^{''}$ is decreasing with $0\leq
J^{''}(0)<+\infty$, by assumption \eqref{drugi moment miary}, we
have {\small
\begin{align*}
&d(t,x)\leq
r_0^{\ast}Be^{\max\Big\{\int_{0}^{t}J^{\prime}\Big(\int_{0}^{t-s+x}\tilde{\lambda}(s,v)r_1(s,v)dv\Big)\tilde{\lambda}(s,t-s+x)ds;\int_{0}^{t}J^{\prime}\Big(\int_{0}^{t-s+x}\tilde{\lambda}(s,v)r_2(s,v)dv\Big)\tilde{\lambda}(s,t-s+x)ds\Big\}}\cdot\\[2ex]
&\cdot\left|\int_{0}^{t}J^{\prime}\left(\int_{0}^{t-s+x}\tilde{\lambda}(s,v)r_1(s,v)dv\right)\tilde{\lambda}(s,t-s+x)ds-\int_{0}^{t}J^{\prime}\left(\int_{0}^{t-s+x}\tilde{\lambda}(s,v)r_2(s,v)dv\right)\tilde{\lambda}(s,t-s+x)ds\right|\\[2ex]
&\leq r_0^{\ast}Be^{\bar{\lambda}T^\ast\max\Big\{
\big|J^{\prime}\big({
\frac{\bar{\lambda}}{\sqrt{\gamma}}\|r_1\|_{\mathbb{L}^{2,\gamma}_{+}}}\big)\big|;
\big|J^{\prime}\big({
\frac{\bar{\lambda}}{\sqrt{\gamma}}\|r_2\|_{\mathbb{L}^{2,\gamma}_{+}}}\big)\big|\Big\}}\cdot
J^{''}(0)\bar{\lambda}^2\int_{0}^{t}\int_{0}^{t-s+x}\mid
r_1(s,v)-r_2(s,v)\mid dv ds\\[2ex]
=& K \int_{0}^{t}\int_{0}^{t-s+x}d(s,v)dv ds, \qquad (t,x)\in
[0,T^\ast]\times[0,+\infty).
\end{align*}}
It follows from Proposition \ref{prop Gronwall} that $r_1=r_2$ on
$[0,T^\ast]\times[0,+\infty)$. \hfill$\square$

\begin{rem}
It follows from Theorem \ref{tw o jedynosci} that under assumptions
of Theorem \ref{tw o istnieniu} $(b)$ the solution is unique in
$\mathbb{H}^{+}_{\gamma}$.
\end{rem}

\section{Explosions}\label{section Explosions}
In this section we show that \eqref{ogr na Jprim} implies that there
is no non-exploding solution of $\eqref{rownanie operatorwe na r}$
in $L^{2,\gamma}_{+}$ on any finite interval $[0,T^\ast]$.

\begin{prop}\label{prop o eksplozjach}
Assume that $J^{\prime}$ satisfies \eqref{ogr na Jprim}. Then for
arbitrary $\kappa \in (0,1)$, there exists a positive constant $K$
such that if
\begin{gather}\label{oddzielenie r_0}
r_0(x)> K, \quad\forall x\in[0,T^\ast],
\end{gather}
then
\begin{gather*}
P(\bar{h}\notin\mathbb{L}^{2,\gamma}_{+})\geq \kappa.
\end{gather*}
\end{prop}
{\bf Proof:} In this proof we use Musiela as well as standard
parametrization. The condition \eqref{oddzielenie r_0} can be
written as $f(0,T)>K$ for $T\in[0,T^\ast]$ and by Theorem 3.4 in
\cite{BarZab} it follows that there is no $f(t,T), 0\leq t\leq
T^\ast, 0\leq T\leq T^\ast$ solving equation \eqref{rownanie
operatorowe na f} which is bounded with
probability grater or equal than $\kappa$.\\
Now assume to the contrary that
$P(\bar{h}\in\mathbb{L}^{2,\gamma}_{+})>1-\kappa$. Due to the
implication
\begin{gather*}
\bar{h}=\mathcal{K}\bar h, \quad
\bar{h}\in\mathbb{L}^{2,\gamma}_{+}\Longrightarrow
\sup_{t\in[0,T^\ast], x\geq0}\bar{h}(t,x)<+\infty,
\end{gather*}
we see that then $\bar{h}$ is bounded with probability grater than
$1-\kappa$. That is a contradiction.\hfill$\square$

\begin{tw}\label{tw o eksplozjach}
Under the assumptions of Proposition \ref{prop o eksplozjach} there
is no solution to the equation \eqref{rownanie operatorwe na r}
taking values in $\mathbb{L}^{2,\gamma}_{+}$ with probability one.
\end{tw}
{\bf Proof:} Assume that $\bar{r}$ is a solution of \eqref{rownanie
operatorwe na r} taking values in $\mathbb{L}^{2,\gamma}_{+}$. Then
$0\leq \bar{r}$ and due to the monotonicity of the operator
$\mathcal{K}$ we see that
\begin{gather*}
h_n(t,x)\leq \bar{r}(t,x), \quad 0\leq t\leq T^\ast, x\geq0, \quad
\forall n=1,2,... \ .
\end{gather*}
Passing to the limit we obtain $\bar{h}\leq\bar{r}$. Thus if
$P(\bar{r}\in\mathbb{L}^{2,\gamma}_{+})=1$ then
$P(\bar{h}\in\mathbb{L}^{2,\gamma}_{+})=1$ which is a contradiction
in view of Proposition \ref{prop o eksplozjach}. \hfill$\square$

\begin{rem}
Due to the inclusion
$\mathbb{H}^{+}_{\gamma}\subseteq\mathbb{L}^{2,\gamma}_{+}$ it
follows that the condition \eqref{ogr na Jprim} implies that there
is no non-exploding solution of \eqref{rownanie operatorwe na r} in
the space $H^{1,\gamma}_{+}$.
\end{rem}

\begin{cor}
One can formulate Theorem \ref{tw o eksplozjach} for other classes
of functions. Below we specify some examples with short
explanations.
\begin{enumerate}[a)]
\item \begin{gather*}
\left\{h: \int_{0}^{+\infty}\mid h(x)\mid^2 w(x)dx<+\infty\right\},
\end{gather*}
where $w:\mathbb{R}_{+}\longrightarrow\mathbb{R}_{+}$ is such that
$\int_{0}^{+\infty}\frac{1}{w(x)}dx<+\infty$. If $h$ belongs to this
space then it is integrable on $(0,+\infty)$. Thus if $r$ is a
non-exploding solution taking values in this space then  must be
bounded which is a contradiction.
\item
\begin{gather*}
\left\{h: \ \mid h(0)\mid + \int_{0}^{+\infty}\mid h^{\prime}(x)\mid^2
w(x) dx<+\infty\right\},
\end{gather*}
where $w:\mathbb{R}_{+}\longrightarrow\mathbb{R}_{+}$ is bounded
from below by $c>0$. If $r$ is a non-exploding solution taking
values in this space then it is locally bounded. Indeed, in view of
the estimation
\begin{align*}
\mid r(t,x)\mid&=r(t,0)+\int_{0}^{x}r^{\prime}(t,v)dv\leq
r(t,0)+\frac{1}{\sqrt{c}}\int_{0}^{x}r^{\prime}(t,v)\sqrt{w(v)}dv\\[2ex]
&\leq r(t,0)+\sqrt{\frac{x}{c}}\sqrt{\int_{0}^{x}\mid
r^{\prime}(t,v)\mid^2w(v)dv}\\[2ex]
&\leq\left(1+\sqrt{\frac{x}{c}}\right)\left(r(t,0)+\sqrt{\int_{0}^{+\infty}\mid
r^{\prime}(t,v)\mid^2 w(v)dv}\right)\\[2ex]
\end{align*}
we have
\begin{gather*}
\sup_{t\in[0,T^\ast], x\in[0,y]}\mid r(t,x)\mid\leq
\left(1+\sqrt{\frac{y}{c}}\right)\sup_{t\in[0,T^\ast]}\left(r(t,0)+\sqrt{\int_{0}^{+\infty}\mid
r^{\prime}(t,v)\mid^2 w(v)dv}\right).
\end{gather*}
Hence $r$ is bounded on each set of the form $[0,T^\ast]\times[0,y]$
and we can proceed as in the proofs of Proposition \ref{prop o
eksplozjach} and Theorem \ref{tw o eksplozjach}.
\end{enumerate}

\end{cor}

\section{Existence results and the L\'evy measure of the noise}\label{Existence and characteristics of the noise}
In this section we gather conditions expressed in terms of the
parameters of the noise which imply \eqref{logarytmiczny wzrost
J^prim} or \eqref{ogr na Jprim}. The proofs can be found in
\cite{BarZab}. To see explicit examples we refer the reader to
\cite{BarZab}.

The first result states that the necessary condition for existence
is that the noise does not contain Wiener part and does not admit
negative jumps, i.e. $q=0$ and $J_1\equiv0$ in \eqref{J
decomposition}.

\begin{tw}\label{tw o warunkach koniecznych}
If the Laplace exponent $J$ of $L$ is such that $q>0$ or
$\nu\{(-\frac{1}{\bar{\lambda}},0)\}>0$ then \eqref{ogr na Jprim}
holds.
\end{tw}

In the following results we assume that $q=0$ and $J_1\equiv0$. It
turns out that then the crucial point is the behavior of the first
derivative of $J_2$ near zero. To formulate the result recall the
concept of slowly varying functions. A positive function $M$ {\it
varies slowly at $0$} if for any fixed $x>0$
\begin{gather*}
\frac{M(tx)}{M(t)}\longrightarrow 1, \qquad \text{as} \
t\longrightarrow0.
\end{gather*}
Typical examples are constants or, for arbitrary $\gamma$ and small
positive $t$, functions
$$
M(t) = \left(\ln {\frac {1}{t} }\right)^{\gamma}\,\,\,.
$$
It turns out that a useful criteria can be formulated in terms of
the behavior near zero of the function
\begin{gather*}
U_{\nu}(x):=\int_{0}^{x}y^2\nu(dy),\qquad x\geq0.
\end{gather*}
Below the notation $f(x)\sim g(x)$ stands for two functions
satisfying
\begin{gather*}
\frac{f(x)}{g(x)}\longrightarrow 1, \qquad \text{as} \
x\longrightarrow0.
\end{gather*}

\begin{tw}\label{tw glowne Tauber}
Assume that for some $\rho\in(0,+\infty)$,
\begin{gather}\label{warunek na zachowanie U w zerze}
U_{\nu}(x)\sim  x^{\rho} \cdot M(x), \qquad as \ x\rightarrow0,
\qquad
\end{gather}
where $M$ is a slowly varying function at $0$.
\begin{enumerate}[i)]
\item If $\rho>1$ then \eqref{logarytmiczny wzrost J^prim} holds.
\item If $\rho<1$, then \eqref{ogr na Jprim} holds.
\item If $\rho=1$, the measure $\nu$ has a density and
\begin{equation}\label{warunek na L}
M(x) \longrightarrow 0 \quad \text{as} \ x\rightarrow 0, \quad and
\quad \int_{0}^{1}\frac{M(x)}{x}\ dx=+\infty,
\end{equation}
then \eqref{logarytmiczny wzrost J^prim} holds.
\end{enumerate}
\end{tw}

\noindent As the next proposition shows, the condition
\eqref{logarytmiczny wzrost J^prim} is satisfied for subordinators
with drifts. This is the special case when the function $J^{\prime}$
is bounded, and thus \eqref{logarytmiczny wzrost J^prim} obviously
holds.

\begin{prop}\label{tw o subordynatorze}
If the process $L$ is a sum of a subordinator and a linear function
then \eqref{logarytmiczny wzrost J^prim} holds. In particular if $L$
is a compound Poisson process with a drift and positive jumps only
then \eqref{logarytmiczny wzrost J^prim} holds.
\end{prop}


\begin{thebibliography}{99}

\bibitem{BarZab}
Barski M., Zabczyk J.: "Forward rate models with linear volatility",
(2010) to appear in {\it Finance and Stochastics},

\bibitem{bjork}
Bj\"ork, Th., Di Masi, G., Kabanov, Y., Runggaldier, W.: "Towards a
general theory of bond markets", (1997), {\it Finance and Stochastics} {\bf 1},
141-174,

\bibitem{FilTap}
Filipovi\'c D., Tappe S.: "Existence of L\'evy term structure
models", (2008), {\it Finance and Stochastic}, 12, 83-115,


\bibitem{BGM} Brace, A., Gatarek, D., Musiela, M.:"The market model of
    interest rate dynamics", (1997), {\it Mathematical Finance}, 7,
    127-147,

\bibitem{Eberlein} Eberlein, E., Raible, S.:"Term structure models
    driven by general L\'evy processes", (1999), {\it Math. Finance}, 9,
    31-53,

\bibitem{EbJacRaib}
Eberlein E., Jacod J., Raible S. "L\'evy term structure models:
No-arbitrage and completeness", (2005), {\it Finance and
Stochastics}, 9, p.67-88,


\bibitem{Filipovic} Filpovi\'c, D.: "Term-Structure Models: A Graduate Course", (2009), Springer-Verlag,

\bibitem{FilipovicTappeTeichmann}
Filpovi\'c, D., Tappe, S., Teichmann, J. : Term structure models
driven by Wiener process and Poisson measures: Existence and
positivity, (2010), {\it SIAM Journal on Financial Mathematics},
Vol.1, 523-554,


\bibitem{HJM}
Heath, D., Jarrow, R., Morton, A.: "Bond pricing and the term
structure of interest rates: a new methodology for contingent claim
valuation",\ (1992), \ {\it Econometrica}, 60, p.77-105,

\bibitem{JakubowskiZabczyk} Jakubowski, J., Zabczyk J.: "Exponential
moments for HJM models with jumps", (2007), {\it Finance and
Stochastics}, 11, 429-445,

\bibitem{Marinelli}
Marinelli, C.: "Local well-posedness of Musiela's SPDE with L\'evy
noise", (2010), {\it Mathematical Finance} {\bf 20}, 341-363,

\bibitem{Morton}
Morton, A.: "Arbitrage and martingales",\ (1989),\ Dissertation,
Cornell University,

\bibitem{PeszatZabczyk}
Peszat, Sz., Zabczyk J.: "Stochastic partial differential equations
with L\'evy noise",(2007), Cambridge University Press,

\bibitem{Rusinek} Rusinek, A.: "Invariant measures for forward rate HJM
    model with L\`evy noise", (2006), Preprint IMPAN 669,
    http://www.impan.pl/Preprints/p669.pdf,


\bibitem{Sato} Sato, K.I.: "L\'evy Processes and Infinite Divisible
    Distributions", (1999), Cambridge University Press.


\end{thebibliography}
\end{document}